\documentclass{amsart}
\usepackage{amsmath,amssymb,mathrsfs}
\usepackage[all]{xy}
 \usepackage[dvips]{graphicx}

\usepackage{comment}

\newtheorem{theo}{Theorem}[section]

\newtheorem{defi}{Definition}[section]

\newcommand{\mX}{{\mathcal X}}
\newcommand{\mL}{{\mathcal L}}
\newcommand{\mcX}{{\overline{\mathcal X}}}
\newcommand{\mcL}{{\overline{\mathcal L}}}
\newcommand{\Td}{{\rm Td}}
\newcommand{\mC}{{\mathbb C}}
\newcommand{\Xinf}{{\mcX_{\infty}}}
\newcommand{\Linf}{{\mcL|_{\mcX_{\infty}}}}

\begin{document}

\title{Cohomology Formula for Obstructions to Asymptotic Chow Semistability}
\author{Yuta Suzuki}
\address{Graduate School of Mathematical Sciences, The University of Tokyo, 3-8-1 Komaba Meguro-ku Tokyo 153-8914, Japan}
\email{ysuzuki@ms.u-tokyo.ac.jp}

\begin{abstract}
Odaka \cite{odaka} and Wang \cite{wang} proved the intersection formula for the Donaldson-Futaki invariant. In this paper, we generalize this result for the higher Futaki invariants, which are obstructions to asymptotic Chow semistability.
\end{abstract}



\maketitle

\section{Introduction}
Let $M$ be a complex manifold of dimension $n$, and $L$ be an ample line bundle over $M$. The existence of constant scalar curvature K\"ahler metrics (cscK metrics for short) in $c_1(L)$ is an important problem in K\"ahler geometry. The conjecture of Yau, Tian and Donaldson asserts that the existence of cscK metrics in $c_1(L)$ is equivalent to a certain GIT stability of the polarized variety $(M,L)$. In fact, Chen, Donaldson and Sun \cite{cds}, Tian \cite{tian} gave the solution of this problem when the polarization is the anti-canonical bundle of $M$. The relevant GIT condition in this case is called K-stability. On the other hand, K-stability is not the only stability notion related to the existence of cscK metrics, and  asymptotic Chow stability is one of stability notions. Donaldson \cite{donaldson01} proved that if the automorphism group ${\rm Aut}(M,L)$ is discrete, then the existence of cscK metrics implies the asymptotic Chow semistability of $(M,L)$. Mabuchi \cite{mabuchi05} extended the result for the case ${\rm Aut}(M,L)$ is not discrete. More precisely, he proved that there exists an obstruction to asymptotic Chow semistability and if the obstruction vanishes, the existence of cscK metrics in $c_1(L)$ implies the asymptotic Chow polystability of $(M,L)$.
After that, the obstruction was reformulated by Futaki \cite{futaki04} as an integral invariants of $(M,L)$ . They are called higher Futaki invariants.

Della Vedova and Zuddas \cite{dellavedovazuddas} showed that the  higher Futaki invariants are closely related to the Chow weight. Let $(\mX,\mL)$ be a test configuration of $(M,L)$. Then the central fiber $(\mX_0, \mL|_{\mX_0})$ is a polarized scheme endowed with a $\mathbb C^*$-action. For an integer $k$, let $\chi(\mX_0, k\mL|_{\mX_0})$ be the Euler-Poincar\'e characteristic of $(\mX_0, k\mL|_{\mX_0})$ and $w(\mX_0, k\mL|_{\mX_0})$ be the weight of the $\mathbb{C}^*$-action on the one-dimensional vector space \\$\otimes_{i=0}^{n}(\wedge^{\rm{max}}H^i(M,L)^{(-1)^i})$. Then the Chow weight ${\rm Chow}(\mX_0, k\mL|_{\mX_0})$ is defined by using $\chi(\mX_0, k\mL|_{\mX_0})$ and $w(\mX_0, k\mL|_{\mX_0})$. We are interested in the asymptotic behavior of ${\rm Chow}(\mX_0, k\mL|_{\mX_0})$ when $k$ grows. Della Vedova and Zuddas showed that if $(\mX_0, k\mL|_{\mX_0})$ is smooth, the higher Futaki invariants are equal to the coefficients of polynomial expansion of ${\rm Chow}(\mX_0, k\mL|_{\mX_0})$ with respect to $k$. The leading coefficient is called the Donaldson-Futaki invariant.

On the other hand, Odaka \cite{odaka} and Wang \cite{wang} showed that there exists another formula of the Donaldson-Futaki invariant. A test configuration $(\mX,\mL)$ has a natural compactification $(\mcX,\mcL)$. Then the Donaldson-Futaki invariant of $(\mX,\mL)$ is expressed by intersection numbers of $(M,L)$ and $(\mX,\mL)$. Theorem \ref{main}, which is our main result, is a generalization of this result.
\begin{theo}\label{main}
Let $(M,L)$ be a polarized variety of complex dimension $n$. Let $(\mathcal X,\mathcal{L})$ be a test configuration of $(M,L)$. If $\mathcal X$ is smooth, then the $\ell$-th Futaki invariant $F_\ell(\mathcal X,\mathcal{L})$ can be computed by the following formula for all $\ell$:
\begin{align*}
F_\ell(\mathcal X,\mathcal L) = \frac{1}{n!(n-\ell+1)!(L^n)^2} \bigg[ &(n+1)(L^n)\big\{(c_1(\overline{\mathcal L} )^{n+1-\ell} {\rm Td}_\ell(\overline{\mathcal X}))\\
 - &(c_1(L)^{n+1-\ell}{\rm Td}_{\ell-1}(M)) \big\} \\
  - &(n-\ell+1)(\overline{\mathcal L }^{n+1})(c_1(L)^{n-\ell}{\rm Td}_\ell(M))     \bigg ],
  \end{align*}
  where $(\overline{\mathcal X},\overline{\mathcal L})$ is the natural compactification of $(\mX,\mL)$.
  \end{theo}
 
 The notation $(\overline{\mathcal L }^{n+1})$ means the intersection number $\mcL\dots\mcL$ in $\mcX$ and $(L^n)$ means $L\dots L$ in $M$, and so on. This theorem allows us to compute $F_\ell(\mathcal X,\mathcal L)$ in terms of characteristic classes of $(M,L)$ and $(\mathcal X,\mathcal L)$.
 
 This paper is organized as follows. In section \ref{background} we recall some definitions and theorems we mentioned above. The proof of Theorem \ref{main} is given in section \ref{proof}. Then we give the compactification of test configurations. In section \ref{localization}, we give the localization formula of Theorem \ref{main}. This is given by Futaki \cite{futaki88}. The localization formula gives an alternative proof for Theorem 1.1 at least for the product configurations. \\
 \noindent {\bf Acknowledgment}
 I am very greatful to my supervisor Akito Futaki for his valuable comments and constant support.

 \section{Background}\label{background}
 \subsection{Asymptotic Chow stability}
 
 First, we recall the definition of asymptotic Chow stability. Let $X \subset \mathbb{P}(V)$ be an $n$-dimensional subvariety of degree $d$. The set of $n+1$ hyperplanes that have a common intersection with $X$
 \begin{align*}
 \{(H_1,\dots,H_{n+1}) \in \mathbb P (V^*)\times\mathbb P (V^*)\times\dots\times\mathbb P (V^*) | H_1 \cap \dots \cap H_{n+1} \cap X \ne \emptyset \}
 \end{align*}
 is an irreducible hypersurface in  $(\mathbb P (V^*))^{n+1}$. The defining polynomial $P_X$ is homogeneous of multi-degree $d=$ deg($X$). The polynomial $P_X \in$ Sym$^d(V)^{\otimes{n+1}}$ is called the Chow form.
 
 Let $(M,L)$ be a polarized variety of dimension $n$. Let $M_k\subset \mathbb{P}(H^0(M,L^k)^*)=\mathbb{P}(V_k)$ be the image of the Kodaira embedding of $M$. For $M_k$ we obtain the Chow form $P_{M_k}\in$ Sym$^{d_k}(V_k)^{\otimes{n+1}}$.
 Consider the ${\rm SL}(V_k)$-action on Sym$^{d_k}(V_k)^{\otimes{n+1}}$.
 \\
 \begin{defi}[Chow stability]
  \begin{enumerate}
\item The variety $M$ is Chow polystable with respect to $L^k$ if the ${\rm SL}(V_k)$-orbit of  $P_{M_k}$ in {\rm Sym}$^{d_k}(V_k)^{\otimes{n+1}}$ is closed.
\item The variety $M$ is Chow stable with respect to $L^k$ if $M$ is polystable and the stabilizer at $P_{M_k}$ is finite.
\item The variety $M$ is Chow semistable with respect to $L^k$ if the closure of the ${\rm SL}(V_k)$-orbit of  $P_{M_k}$ in {\rm Sym}$^{d_k}(V)^{\otimes{n+1}}$ does not contain the origin ${\bf o}\in {\rm Sym}^{d_k}(V)^{\otimes{n+1}}$ .
  \end{enumerate}
 \end{defi}

\begin{defi}[asymptotic Chow stability]
The variety $M$ is asymptotically Chow polystable (respectively stable or semistable) with respect to $L$ if there is a $k_0>0$ such that for any $k>k_0$, $M$ is Chow polystable (respectively stable or semistable) with respect to $L^k$.
\end{defi}

\subsection{The relationship to cscK metrics}
There is a relation between GIT stabilities and canonical metrics, called the Kobayashi-Hitchin correspondence. The well-known conjecture of Yau, Tian and Donaldson asserts that the polarized manifold $(M,L)$ is ``K-polystable'' if and only if constant scalar curvature K\"ahler metrics (cscK metrics for short) exist in $c_1(L)$. K-stability is a little different notion from Chow stability, but asymptotic Chow semistability implies K-semistability. Chow stability also has a relationship to cscK metric, as explained below.

 We denote by Aut$(M)$ the group of automorphisms of $M$ and by Aut$(L)$ the group of bundle automorphisms of $L$. Let Aut$(M,L)$ be the subgroup of automorphism group Aut$(L)$ of $L$  consisting of all automorphisms of L commuting with the $\mathbb{C}^*$-action on the fiber. Such automorphisms of $L$ descend to automorphisms of $M$. So we can consider Aut$(M,L)$ as a subgroup of Aut$(M)$. Donaldson proved the following theorem.
 
 \begin{theo}[Donaldson\cite{donaldson01}]
 Let $(M,L)$ be a polarized manifold. Assume that Aut$(M,L)$ is discrete. If $M$ admits cscK metrics in $c_1(L)$, then $(M,L)$ is asymptotically Chow stable.
 \end{theo}
The assumption for Aut$(M,L)$ means the finiteness of the stabilizer. This theorem gives us a differential geometric criterion of asymptotic Chow stability. Note that we can not omit the assumption. The following example is known.
 \begin{theo}[Ono-Sano-Yotsutani \cite{onosanoyotsutani}, Nill-Paffenholtz\cite{nillpaffen}] There is a toric Fano $7$-manifold which admits K\"ahler-Einstein metrics (so cscK metrics) but not asymptotically Chow semistable.
\end{theo}
We will explain the example of Ono-Sano-Yotsutani in Section 4.

In the case when Aut$(M,L)$ is not discrete, Mabuchi extended the theorem of Donaldson.

\begin{theo}[Mabuchi\cite{mabuchi04}, \cite{mabuchi05}] If Aut$(M,L)$ is not discrete, there exists an obstruction to asymptotic Chow semistability. If the obstruction vanishes, then the existence of cscK metrics in $c_1(L)$ implies the asymptotic Chow polystability of $(M,L)$.
 
\end{theo}
This obstruction was reformulated by Futaki \cite{futaki04}. We will explain that in the next section.
\subsection{Higher Futaki invariants}
Here we recall the definition of the Futaki invariant. First let $\mathfrak h(M)$ be the complex Lie algebra of Aut$(M)$ which consists of holomorphic vector fields over $M$. When  we consider Aut$(M,L)$ as the Lie subgroup of Aut$(M)$, its Lie algebra $\mathfrak h_0(M)$ is a Lie subalgebra of $\mathfrak h(M)$. Second we fix a K\"ahler form $\omega$ representing $c_1(L)$. Then for any $X\in\mathfrak h_0$ there exists a complex valued function $u_X$ such that 
\begin{align}
i(X)\omega=-\overline\partial u_X,\\ \int_M u_X\ \omega^n = 0.\label{normalization}
\end{align}
The function $u_X$ is called the Hamiltonian function of $X$. The existence of such $u_X$ is well known, see \cite{lebrunsimanca93}. Let $\nabla$ be the Chern connection of the K\"ahler metric associated to $\omega$, and $\Theta$ be the curvature of $\nabla$. Put $L(X) =$ $\nabla{_X}-\mathcal L{_X}=\nabla{X}$ where $\mathcal L{_X}$ is the Lie derivative. The operator $L(X)$ defines a smooth section of ${\rm End}(T^{1,0}M)$. Let $\phi$ be a ${\rm GL}(n,\mathbb C)$-invariant polynomial of degree $\ell$ on $\mathfrak{gl}(n,\mathbb C)$. We define $\mathcal F_{\phi}\colon \mathfrak h_0(M) \rightarrow \mathbb C$ by 
\begin{equation*}
\mathcal F_{\phi}(X) = (n-\ell+1)\int_M\phi(\Theta)\wedge u_X\omega^{n-\ell} + \int_M\phi(L(X) + \Theta)\wedge\omega^{n-\ell+1}.
\end{equation*}
Then $\mathcal F_{\phi}$ does not depend on the choice of $\omega$ and depends only on the K\"ahler class $c_1(L)$ (see \cite{futaki04}). Let Td$_{\ell}$ be the $\ell$-th Todd polynomial, which is $GL(n,\mathbb C)$-invariant polynomial of degree $\ell$ on $\mathfrak{gl}(n,\mathbb C)$. Then $\mathcal F_{\rm{Td}_{\ell}}$ is an obstruction to the asymptotically Chow semistability.
\begin{theo}[Futaki \cite{futaki04}] If $(M,L)$ is asymptotically Chow semistable, then 
\begin{equation*}
\mathcal{F}_{\rm{Td}_{\ell}}(X) = 0
\end{equation*}
holds for any $1 \le \ell \le n$ and $X$ in a maximal reductive subalgebra $\mathfrak h_r$ of $\mathfrak h_0(M)$. The vanishing of all invariants $\mathcal F_{\rm{Td}_{\ell}}$ is equivalent to the vanishing of Mabuchi's obstruction.
\end{theo}
The invariant $\mathcal F_{\rm{Td}_{\ell}}$ is called the $\ell$-th Futaki invariant. In paticular first Futaki invariant $\mathcal F_{\rm{Td}_{1}}$ is the same as classical one up to a constant factor.

\subsection{Chow weight}
There is another interpretation of $\mathcal{F}_{\rm{Td}_{\ell}}$ by Della Vedova and Zuddas \cite{dellavedovazuddas}. Given a one-parameter subgroup $\rho \colon \mathbb C ^*\rightarrow{\rm Aut}(M,L)$ with a lifting action on $L$, let $X  \in \mathfrak h_{0}(M)={\rm Lie(Aut}(M,L)\rm{)}$ be the generator of $\rho$. We denote by $w(M,L)$ the weight of the $\mathbb C^*$-action induced on $\otimes_{i=0}^{n}(\wedge^{\rm{max}}H^i(M,L)^{(-1)^i})$, and by $\chi(M,L)$ the Euler-Poincar\'e characteristic $\sum_{i=0}^n(-1)^i\dim H^i(M,L)$. For sufficiently large $k$, we may assume $H^i(M,L^k) = 0$ for $i>0$ by the Kodaira vanishing theorem. We have        polynomial expansions with respect to $k$: 
\begin{align}
\chi(M,L^k) &= a_0(M,L)k^n + a_1(M,L)k^{n-1} + \dots + a_n(M,L), \\
w(M,L^k) &= b_0(M,L)k^{n+1} + b_1(M,L)k^n + \dots + b_{n+1}(M,L).
\end{align}

\begin{defi}The Chow weight of this action is defined by
 \begin{align*}
{\rm Chow}(M,L^k) = \frac{w(M,L^k)}{k \chi (M,L^k)} - \frac{b_0(M,L)}{a_0(M,L)}.
 \end{align*} 
\end{defi}
We can show
\begin{align*}
{\rm Chow}(M,L^k) &= \frac{b_{n+1}(M,L)}{\chi(M,L^k)} \\
                  &+ \frac{a_0(M,L)}{k \chi(M,L^k)}\sum_{\ell=1}^{n}\frac{a_0(M,L)b_{\ell}(M,L)-b_0(M,L)a_{\ell}(M,L)}{a_0(M,L)^2} k^{n-\ell+1}.
\end{align*}
The first term is known to vanish in the smooth case. Therefore we define $F_{\ell}(M,L)$ by
\begin{align}
F_{\ell}(M,L) = \frac{a_0(M,L)b_{\ell}(M,L)-b_0(M,L)a_{\ell}(M,L)}{a_0(M,L)^2}.
\end{align}
\begin{theo}[Della Vedova-Zuddas \cite{dellavedovazuddas}]
\begin{align}
F_{\ell}(M,L) = \frac{1}{Vol(M,L)}\mathcal F_{{\rm Td}_{\ell}}(X).
\end{align}
\end{theo}
Paul and Tian showed that the first Futaki invariant $F_1$ can be considered as the Mumford weight of the CM-line $\lambda_{\rm Chow}$ on the Hilbert scheme. Della Vedova and Zuddas showed that the $\ell$-th Futaki invariant is also the weight of some line $\lambda_{{\rm Chow},\ell}$, see \cite{dellavedovazuddas} and \cite{futaki12} for detail.

\subsection{Intersection formula of Donaldson-Futaki invariant}
So far we have considered product configurations. Now we consider general test configurations. First we recall the definition. Note that Tian originally assumed that a central fiber of a test configuration $\mathcal X_0$ is normal but Donaldson defined a test configuration as a general scheme. After that Li and Xu proved we can assume its normality by using the minimal model program (see \cite{lixu}).

\begin{defi}[test configuration]\label{testconfiguration}
Let $(M,L)$ be a polarized variety. A test configuration of $(M,L)$ consists of the following data:
\begin{enumerate}
\item a scheme $\mathcal X$ with $\mathbb C^*$-action;
\item a $\mathbb C^*$-equivariant relative ample line bundle $\mathcal L \rightarrow \mathcal X$;
\item a flat $\mathbb C^*$-equivariant morphism $\pi \colon (\mathcal X,\mathcal L) \rightarrow \mathbb C$ such that\\ we have $(\mathcal X_1,\mathcal L|_{\mathcal X_1}) = (M,L)$ where $\mathcal X_1 = \pi^{-1}(1)$.
\end{enumerate}
\end{defi}

When $(M,L)$ has a $\mathbb C^*$-action, we can make a test configuration $\mathcal X = M\times\mathbb C$ with a diagonal $\mathbb C^*$-action, called a product configuration.  Clearly the central fiber $(\mathcal X_0, \mathcal L|_{\mathcal X_0})$ is a polarized scheme endowed with a $\mathbb C^*$-action .

\begin{defi}[Donaldson-Futaki invariant] Let $(\mathcal X,\mathcal L)$ be a test configuration of $(M,L)$. Then the Donaldson-Futaki invariant DF$(\mathcal X,\mathcal L)$ is defined by $F_1(\mathcal X_0,\mathcal L|_{\mathcal X_0})$.
\end{defi}

The intersection formula for the Donaldson-Futaki invariant DF$(\mathcal X,\mathcal L)$ was proved by Odaka and Wang.

\begin{theo}[Odaka\cite{odaka}, Wang\cite{wang}] If $(\mathcal X,\mathcal L)$ is normal and $\mathbb Q$-Gorenstein, it follows that
\begin{align*}
{\rm DF}(\mathcal X,\mathcal L) = \frac{1}{2(n+1)(L^n)^2}\Big((n+1)(K_{\overline{\mathcal X}/\mathbb P^1}.\overline{\mathcal L}^n)(L^n)-n(\overline {\mathcal L}^{n+1})(K_M.L^{n-1})\Big),
\end{align*}
where $(\overline{\mathcal X},\overline{\mathcal L})$ is the  natural compactification of $(\mathcal X,\mathcal L)$, explained in the proof of our Theorem\ref{main}.
\end{theo}
Theorem \ref{main} is a generalization of this result.

\section{Proof of main result}\label{proof}
First we recall the compactification of $(\mathcal X,\mathcal L)$ in \cite{wang}. Let $[z_0:z_1]$ be the homogeneous coordinates of $\mathbb C \mathbb P^1$. Let $0=[1:0]$, $\infty = [0:1]$, $\Delta_0$, $\Delta_{\infty}$ be the coordinate neighborhood of $0$, $\infty$ respectively and $\mu = z_1/z_0$ be the local coordinate in $\Delta_0$. The transition function $g\colon \Delta_0\backslash 0 \rightarrow \Delta_{\infty}\backslash 0$ is given by $g(\mu) = 1/\mu$.

We define a $\mathbb C^*$-action on $\mathbb C \mathbb P^1 \times M$ by
\begin{equation}
t\cdot([z_0:z_1], p) = ([z_0: t z_1],p)
\end{equation}
for $t \in \mathbb C^*$. In $\Delta_0 \times M$, this action is given by $t.(\mu, p) = (t\mu, p)$.
Let $x \in \mathcal X \backslash \mathcal X_0$ and $\pi(x)=\mu$. Then
\begin{align*}\label{yamamoto}
f\colon \mathcal X \backslash \mathcal X_0 \rightarrow \Delta_0\backslash \{0\} \times M \\
x \mapsto (\mu, \rho(\mu)^{-1}x)
\end{align*}
is isomorphic where $\rho \colon \mathbb C^* \rightarrow {\rm Aut}(\mathcal X, \mathcal L)$ is the $\mathbb C^*$-action and we use $\rho(\mu)^{-1}x \in \mathcal X_1 = M$. Moreover this map is $\mathbb C^*$-equivariant, that is, the following diagram commutes.
\begin{equation*}
\begin{xy}
(0,40) *{\mathcal X\backslash \mathcal X_0}="",
 (40,40) *{\Delta_0\backslash \{0\} \times M}="",
(0,38) *{}="",
(40,38) *{}="",
 (0,30) *{x}="A", (40,30)*{(\mu, \rho(\mu)^{-1}x)}="B",
(-2,0) *{\rho(t)x}="C", (40,0)*{(t \mu, \rho(\mu)^{-1}x)}="D",

\ar^f @{|->} (3,30);(28,30)
\ar^{\times t} @{|->} (0,28);(0,3)
\ar^{f} @{|->} (3,0);(28,0)
\ar^{\times t} @{|->} (40,27);(40,3)
\end{xy}
\end{equation*}

So we can define

\begin{equation*}
\overline{\mathcal X} =\mathcal X_0 \cup (\mathcal X \backslash \mathcal X_0) \bigcup_{(g\times id_M) \circ f} \Delta_{\infty} \times M 
\end{equation*}
with the $\mathbb C^*$-action. By $\mathcal L|_{\mathcal X \backslash \mathcal X_0} = \Delta_0 \backslash 0 \times L$, we also have the $\mathbb C^*$-equivariant line bundle $\overline{\mathcal L} \rightarrow \overline{\mathcal X}$. Finally we get the $\mathbb C^*$-equivariant morphism $\pi \colon (\overline{\mathcal X}, \overline{\mathcal L}) \rightarrow \mathbb C \mathbb P^1$.
Note that the $\mathbb C^*$-action on $\overline{\mathcal X}_{\infty}=\pi^{-1}(\{\infty\})$ and $\overline{\mathcal L}|_{\mathcal X_{\infty}}$ is trivial.

We can take an integer $N$ such that $\overline{\mathcal L} + N\pi^*(\mathcal O_{\mathbb C\mathbb P^1}(1)) \rightarrow \overline{\mathcal X}$ is ample since $\overline{\mathcal L} \rightarrow \overline{\mathcal X}$ is relatively ample. Put
\begin{equation}\label{N}
\mathcal N = \overline{\mathcal L} + N\pi^*\mathcal O_{\mathbb C\mathbb P^1}(1).
\end{equation}
 Let $\sigma_0$ and $\sigma_{\infty}$ be the sections of $\mathcal O_{\mathbb{CP}^1}(1)$ corresponding to the divisors $[0]$ and $[\infty]$, respectively. Since the action on $\mathbb{CP}^1$ lifts to $\mathcal O_{\mathbb C\mathbb P^1}(-1)$, the weight of $\sigma_0$ and $\sigma_{\infty}$ are $-1, 0$, respectively.
The short exact sequence
\[\xymatrix{
0 \ar[r]
& k\mathcal N(-[\mathcal X_0]) \ar^-{\times \pi^*{\sigma_0}}[r] 
& k \mathcal N\ar[r] 
& k\mathcal N|_{\mathcal X 0} \ar[r] 
&0
}\]
induces the following exact sequence
\[\xymatrix{
0 \ar[r]
& H^0(\overline{\mathcal X},k\mathcal N(-[\mathcal X_ 0])) \ar^-{\times \pi^*{\sigma_0}}[r]
&  H^0(\overline{\mathcal X},k\mathcal N) \ar[r]
& H^0(\mathcal X_0,k\mathcal N|_{\mathcal X_0}) \ar[r]
& 0
}\]
by the Kodaira vanishing theorem. Put $V_1 = H^0(\overline{\mathcal X},k\mathcal N(-[\mathcal X_ 0]))$, 
$V_2 = H^0(\overline{\mathcal X},k\mathcal N)$  and $V_3 = H^0(\mathcal X_0,k\mathcal N|_{\mathcal X_0})$. We denote by $d_i$ the dimension of $V_i$ and by $w_i$ the total weights of the ${\mathbb C}^*$-action on $V_i$. We have
\begin{align}
d_3 &= d_2 - d_1,\label{d1}\\
w_3 &= w_2 - (w_1 - d_1)\label{w1},
\end{align}
since the weight of $\pi^*\sigma_0$ is $-1$. By definition, the $\mathbb C^*$-action on $\pi^* \mathcal O_{\mathbb{CP}^1}(1)|_{\mathcal X_0}$ is trivial, and hence $w_3 =  w(\mX_0,k\mL|_{\mX_0})$. For $\sigma_{\infty}$, we have 
\[\xymatrix{
0 \ar[r]
& H^0(\overline{\mathcal X},k\mathcal N(-[\overline{\mathcal X}_ \infty])) \ar^-{\times \pi^*{\sigma_{\infty}}}[r]
&  H^0(\overline{\mathcal X},k\mathcal N) \ar[r]
& H^0(\overline{\mathcal X}_{\infty},k\mathcal N|_{\overline{\mathcal X}_{\infty}}) \ar[r]
& 0.
}\]
Put $V_4 = H^0(\overline{\mathcal X}_{\infty},k\mathcal N|_{\overline{\mathcal X}_{\infty}})$. Similarly, we have
\begin{align}
d_4 &= d_2 - d_1,\label{d2} \\
w_4 &= w_2 - w_1.\label{w2}
\end{align}
Note that $k\mathcal N|_{\overline{\mathcal X}_{\infty}}= k\overline{\mathcal L}|_{\overline{\mathcal X}_{\infty}} + kN\pi^*\mathcal O _{\mathbb P^1}([\infty]) $ by (\ref{N}). Since the action on $\overline{\mathcal L}|_{\overline{\mathcal X}_{\infty}}$ is trivial, we have $w_4 = -kNd_4$.
It follows that
\begin{align}\label{weight}
w(\mX_0,k\mL|_{\mX_0}) = w_3 
       &=w_2 - (w_1 - d_1) \nonumber \\
       &=d_1 + (w_2 - w_1) \nonumber \\
       &=d_1 + w_4  \nonumber \\ 
       &=d_2 - (kN + 1)d_3 \nonumber \\
       &=\dim H^0(\overline{\mathcal X},k \mathcal  N )-(kN+1)\dim H^0(\mX_0,k\mL|_{\mX_0})
\end{align}
from (\ref{d1}) to  (\ref{w2}). 

Now, we calculate the weight $w(\mX_0,k \mL|_{\mX_0})$ by the Riemann-Roch-Hirzebruch theorem. Note that $\dim H^0(\mX_0, k \mL|_{\mX_0})  =\dim H^0(\mX_1,k \mL|_{\mX_1}) = \dim H^0(M,kL)$  for sufficiently large $k$ by the flatness of $\mcX \rightarrow \mathbb {CP}^1$.
We have
\begin{align}\label{ML}
\begin{aligned}
\dim H^0(M,kL)&=\int_M {\rm ch}(kL){\rm Td}(M)\\
              &=\sum_{\ell = 0}^n \frac{1}{(n-\ell)!}\int_M c_1(L)^{n-\ell}{\rm Td}_{\ell}(M)k^{n-\ell},
\end{aligned}
\end{align}
\begin{align}\label{XN}
\dim H^0(\overline{\mathcal X},k{\mathcal N})&=\int_{\mathcal X} {\rm ch}(k\mathcal N){\rm Td}(\overline{\mathcal X})\nonumber \\
              &=\sum_{\ell = 0}^{n+1} \frac{1}{(n-\ell+1)!}\int_{\mcX}c_1(\mathcal N)^{n-\ell+1}{\rm Td}_{\ell}(\overline{\mathcal X})k^{n-\ell+1}\nonumber \\
              &=\sum_{\ell = 0}^{n+1} \frac{1}{(n-\ell+1)!}\int_{\mcX}c_1(\overline{\mathcal L})^{n-\ell+1}{\rm Td}_{\ell}(\overline{\mathcal X}) k^{n-\ell +1}\nonumber \\
              &+ N\sum_{\ell = 0}^{n} \frac{1}{(n-\ell)!}\int_{\mathcal X_{\infty}} c_1(\overline{\mathcal L})^{n-\ell}{\rm Td}_{\ell}(\overline{\mathcal X})k^{n-\ell+1}.
\end{align}
Here we use $\mathcal N = \overline{\mathcal L} + N\pi^*\mathcal O_{\mathbb C\mathbb P^1}([\infty])$. Substituting (\ref{ML}) and (\ref{XN}) into (\ref{weight}), we get
\begin{align*}
&w(\mX_0,k\mL|_{\mX_0}) = w_3 \\&= \Big[\frac{1}{(n+1)!}c_1(\mcL)^{n+1} + N\frac{1}{n!} \big\{\int_{\mcX_{\infty}}c_1(\mcL)^n-\int_M c_1(L)^n\big\}\Big]k^{n+1}\\
&+ \sum_{\ell=1}^{n}\Big[\frac{1}{(n-\ell+1)!} \{c_1(\overline{\mathcal L})^{n-\ell+1}{\rm Td}_{\ell}(\overline{\mathcal X})-c_1(L)^{n-\ell+1}\Td_{\ell -1}(M) \}\\  &+N\frac{1}{(n-\ell)!}\big\{\int_{\mcX_{\infty}}c_1(\mcL)^{n-\ell}\Td_{\ell}(\mcX) - \int_M c_1(L)^{n-\ell}\Td_{\ell}(M)\big\}\Big] k^{n-\ell +1}\\
&+ \Big[\int_{\mcX}\Td_{n+1}(\mcX)-\int_M\Td_n(M) \Big].
\end{align*}
Note that this polynomial does not depend on $N$. Finally, we obtain
\begin{align}\label{bl1}
b_{\ell} = \frac{1}{(n-\ell+1)!}\Big[c_1(\mcL)^{n-\ell+1}\Td_{\ell}(\mcX)-c_1(L)^{n-\ell+1}\Td_{\ell-1}(M)\Big]
 \end{align}
for $1 \leq \ell \leq n$. 
This implies Theorem~\ref{main}.

\section{Localization and Example}\label{localization}
In this section, we will see that Theorem~\ref{main} is localized to the formula in \cite{futaki88} by the original Bott residue formula. This gives the alternative proof of the result of Della Vedova and Zuddas \cite{dellavedovazuddas}. Finally we give the example of \cite{nillpaffen} calculated in \cite{onosanoyotsutani}.

For the convenience of reader, recall the Bott residue formula \cite{bott}. Let $M$ be a compact complex manifold and $\varphi$ be a ${\rm GL}(n,\mathbb C)$-invariant polynomial of degree $n$ on $\mathfrak{gl}(n,\mathbb C)$. Let $X$ be a holomorphic vector field. Assume that the zero set of X consists of manifolds $\{ Z_{\lambda}\}$. Then $L(X) = \nabla_X-\mathcal L_X$ induces a endomorphism $L^{\nu}(X)$ of the normal bundle $\nu(Z_{\lambda})$. Suppose that $L^{\nu}(X)$ is non-degenerate.
Then it holds
\begin{align}
\varphi(M) = \varphi(\Theta) = \sum_{\lambda \in \Lambda}\int_{Z_{\lambda}}\frac{\varphi(L(X) + \Theta)|_{Z_{\lambda}}}{\det(\frac{\sqrt{-1}}{2\pi}(L^{\nu}(X) + K))},
\end{align}
where $\Theta$ and $K$ is the curvature of tangent bundle $TM$ and  normal bundle $\nu(Z_{\lambda})$ respectively. Note that Bott proved this for arbitrary equivariant vector bundle, not only for tangent bundle $TM$. 

We consider the localization of Theorem \ref{main} by the Bott residue formula. Let $M$ be a Fano manifold and $X$ a holomorphic vector field on $M$. Assume that zero set of $X$ consists of isolated points and $X$ is non-degenerate. Take the canonical lift of $X$ to $-K_M$. Let $(\mX,\mL)$ be the product configuration of the $\mathbb C^*$-action generated by $X$ and $(\mcX,\mcL)$ be the compactification of $(\mX,\mL)$. The set of fixed points of the $\mC^*$-action on $\mcX$ is the union of the whole fiber $\mcX_{\infty}$ and points on the central fiber $\mX_0$. Thus, we have
\begin{align}\label{bl}
\int_{\mcX}c_1(\mcL)^{n-\ell+1}\Td_{\ell}(\mcX) = \int_{\mcX_{\infty}} \frac{c_1(\mcL|_{\mcX_{\infty}})^{n-\ell+1}\Td_{\ell}(\bar{L}(X) + \Theta)}{\det\frac{\sqrt{-1}}{2\pi}(\bar{L}^{\nu}(X) + K)} \nonumber \\+ 
\sum_{{\bf q}:{\rm fixed \  point}} \frac{(c_1^{n-\ell+1}\Td_{\ell})(\bar{L}(X)_{\bf q})}{\det\frac{\sqrt{-1}}{2\pi}(\bar{L}(X)_{\bf q})},
\end{align}
where 
$\bar{L}(X)$ is the endomorphism of tangent bundle $T\mcX$, 
$\bar{L}^{\nu}(X)$ is the induced endomorphism of the normal bundle $\nu(\mcX_{\infty})$, $\Theta$ is curvature of $T\mcX$ and $K$ is the curvature of $\nu(\mcX_{\infty})$. Here we use the fact that $\mcL|_{\Xinf}$ is the anticanonical bundle and the $\mC^*$-action on $\mcL|_{\Xinf}$ is trivial. We consider the first term of (\ref{bl}). we omit the determinant since the codimension of $\Xinf$ is one. From the construction of $\mcX$, $\nu(\Xinf)$ is trivial. Thus, $K=0$ and $T\mcX|_{\Xinf}$ is decomposed to $ T\mathbb P^1|_{\{\infty\}}\oplus   T\Xinf$. Then we have
\begin{equation}\label{linf}
\bar{L}(X) + \Theta = \left(
\begin{array}{c|ccc}
-1     & 0       & \ldots  & 0    \\ \hline
0      &         &         &   \\
\vdots &         &\Theta_M &   \\
0      &         &         &
\end{array} \right),
\end{equation}
where $\Theta_M$ is the curvature of $M$. Since the Todd polynomial is multiplicative, it follows that
\begin{align}\label{loc.todd}
\Td_{\ell}(\bar{L}(X) + \Theta) 
&= \Td_{\ell}(M) + \Td_1(-1) \Td_{\ell-1}(\Theta_M) \nonumber \\
&= \Td_{\ell}(M) + \frac{1}{2}c_1(-1) \Td_{\ell-1}(\Theta_M) \nonumber \\
&= \Td_{\ell}(M) - \frac{\sqrt{-1}}{4\pi}\Td_{\ell-1}(\Theta_M).
\end{align}
Substituting (\ref{linf}) and (\ref{loc.todd}) to (\ref{bl}), we get
\begin{align}\label{cal1}
&\int_{\mcX}c_1(\mcL)^{n-\ell+1}\Td_{\ell}(\mcX)\nonumber \\ &=  \frac{1}{2
}\int_{\Xinf}c_1(\Linf)^{n-\ell+1}\Td_{\ell-1}(\Theta_M)+\sum_{q:{\rm fixed \ point}} \frac{(c_1^{n-\ell+1}\Td_{\ell})(\bar{L}(X))}{\det\frac{\sqrt{-1}}{2\pi}(\bar{L}(X))} \nonumber\\
&=\frac{1}{2}\int_{M}c_1(L)^{n-\ell+1}\Td_{\ell-1}(M) +\sum_{q:{\rm fixed\ point}} \frac{(c_1^{n-\ell+1}\Td_{\ell})(\bar{L}(X))}{\det\frac{\sqrt{-1}}{2\pi}(\bar{L}(X))}.
\end{align}
Similarly we can calculate the second term of (\ref{cal1}). On the central fiber $\mX_0$, we have
\begin{equation}\label{cal2}
\bar{L}(X)  = \left(
\begin{array}{c|ccc}
1      & 0       & \ldots  & 0    \\ \hline
0      &         &         &   \\
\vdots &         &L(X) &   \\
0      &         &         &
\end{array} \right)
\end{equation}
and
\begin{align}\label{cal3}
\Td_{\ell}(\bar{L}(X)) = \Td_{\ell}(L(X)) + \frac{\sqrt{-1}}{4\pi}\Td_{\ell-1}(L(X)).
\end{align}
Substituting (\ref{cal2}) and (\ref{cal3}) to (\ref{cal1}), we obtain
\begin{align}\label{cal4}
&\int_{\mcX}c_1(\mcL)^{n-\ell+1}\Td_{\ell}(\mcX) \nonumber\\&=
\int_{\Xinf}c_1(\Linf)^{n-\ell+1}\Td_{\ell-1}(\Theta)+ \sum_{q:{\rm fixed \ point}} \frac{(c_1^{n-\ell+1}\Td_{\ell})(L(X)_{\bf q})}{\det\frac{\sqrt{-1}}{2\pi}(L(X)_{\bf q})}. \nonumber\\
\end{align}
From (\ref{bl1}), it follows that
\begin{equation}\label{bl2}
b_{\ell} =\frac{1}{(n-\ell+1)!}\sum_{{\bf q}:{\rm fixed\ point}} \frac{(c_1^{n-\ell+1}\Td_{\ell})(L(X)_{\bf q})}{\det\frac{\sqrt{-1}}{2\pi}(L(X)_{\bf q})}.
\end{equation}
This is the localization formula in \cite{onosanoyotsutani}.

Similarly we obtain the localization formula of $b_0$ and $a_{\ell}$ :
\begin{align} 
b_0 &= \frac{(\mcL^{n+1})}{(n+1)!} \nonumber \\
    &= \frac{1}{(n+1)!}\sum_{{\bf q}:{\rm fixed\ point}}\frac{c_1(L(X)_{\bf q})^{n+1}}{\det\frac{\sqrt{-1}}{2\pi}(L(X)_{\bf q})},\label{b_zero}\\
    a_{\ell} &= \frac{1}{(n-\ell)!}\sum_{{\bf q}:{\rm fixed\ point}}\frac{c_1(L(X)_{\bf q})^{n-\ell}\Td_{\ell}(L(X)_{\bf q})}{\det\frac{\sqrt{-1}}{2\pi}(L(X)_{\bf q})}.\nonumber
\end{align}

Next, see the example in \cite{nillpaffen}. We consider the Fano polytope in $\mathbb R^7$ whose vertices are given by the following matrix:
\begin{equation}\label{vertex}
\begin{pmatrix}
1 & 0 & 0 & 0 & 0 & -1 & 0 & 0 & 0 & 0 & 0 & 0\\
0 & 1 & 0 & 0 & -1 & 0 & 0 & 0 & 0 & 0 & 0 & 0\\
0 & 0 & 1 & -1 & 0 & 0 & 0 & 0 & 0 & 0 & 0 & 0\\
0 & 0 & 0 & 0 & 0 & 0 & 1 & 0 & 0 & -1 & 0 & 0\\
0 & 0 & 0 & 0 & 0 & 0 & 0 & 1 & 0 & -1 & 0 & 0\\
0 & 0 & 0 & 0 & 0 & 0 & 0 & 0 & 1 & -1 & 0 & 0\\
0 & 0 & 0 & -1 & -1 & -1 & 0 & 0 & 0 & 2 & 1 & -1
\end{pmatrix}.
\end{equation}
Let $M$ be the $7$-dimensional toric Fano manifold associated with the polytope. Then $M$ is a $\mathbb P^1$-fibration on $(\mathbb P^1)^3\times \mathbb P^3$ and admits K\"ahler-Einstein metrics (see \cite{nillpaffen}).

Note that  $b_0$ in (\ref{b_zero}) coincides with the original Futaki invariant when we consider the canonical lift of $X$ to $-K_M$ (see \cite{futakimorita}). Now it is zero since $M$ is K\"ahler-Einstein. So we just have to calculate $b_{\ell}$ and $a_0$.

Next, define a $\mathbb C^*$-action on $M$. 
Here we consider the following one-parameter subgroup. Let $v_i$ be the $i$-th vertex in (\ref{vertex}). Let Spec$(\mathbb C [X_1,X_2,X_3,Y_1,Y_2,Y_3,Z])$ be the affine toric variety which corresponds to the $7$-dimensional cone generated by $\{v_1, v_2, v_3, v_7, v_8, v_9, v_{11}\}$. Here $X_1, X_2, X_3$ are affine coordinates of $(\mathbb P^1)^3$, $Y_1, Y_2, Y_3$ are affine coordinates of $\mathbb P^3$ and $Z$ is an affine coordinate of the fiber. Then the one-parameter subgroup $\sigma_t$ is defined by
\begin{align*}
\sigma_t\cdot&(X_1,X_2,X_3,Y_1,Y_2,Y_3,Z)\\=&(e^{\alpha_1t}X_1,e^{\alpha_2t}X_2,e^{\alpha_3t}X_3,e^{\beta_1t}Y_1,e^{\beta_2t}Y_2,e^{\beta_3t}Y_3,e^{\gamma t}Z)
.\end{align*}
This one-parameter subgroup is defined over the whole $M$. For a generic \\$\{\alpha_i, \beta_j, \gamma\}_{1\leq i, j \leq3}$, the set of fixed points of $\sigma_t$ consists of the isolated $64$ points (see \cite{onosanoyotsutani}). Let $X$ be the holomorphic vector field generated by $\sigma_t$. Take the lift of X to $-K_M$ as section 2.3. 
Then higher Futaki invariants are calculated in \cite{onosanoyotsutani} using the localization formula (\ref{bl2}). We have

\begin{align*}
b_2 &= \frac{68}{45}(\sum_{\i=1}^{3}\alpha_i - \sum_{\i=1}^{3} \beta_i -2\gamma),\\
b_3 &= \frac{68}{15}(\sum_{\i=1}^{3}\alpha_i - \sum_{\i=1}^{3} \beta_i -2\gamma),\\
b_4 &= \frac{49}{9}(\sum_{\i=1}^{3}\alpha_i - \sum_{\i=1}^{3} \beta_i -2\gamma),\\
b_5 &= \frac{10}{3}(\sum_{\i=1}^{3}\alpha_i - \sum_{\i=1}^{3} \beta_i -2\gamma),\\
b_6 &= \frac{214}{315}(\sum_{\i=1}^{3}\alpha_i -\sum_{\i=1}^{3} \beta_i -2\gamma),\\
b_7 &= \frac{2}{15}(\sum_{\i=1}^{3}\alpha_i - \sum_{\i=1}^{3} \beta_i -2\gamma).
\end{align*}
The similar calculation gives 
\begin{align}
(-K_M)^7 = 13047715.
\end{align}
Finally, we obtain
\begin{align*}
F_2(\mX,\mL) 
&= \frac{ a_0 b_{2} - b_0 a_{2}}{a_0 ^2} \\
&= \frac{7616}{13047715}(\sum_{\i=1}^{3}\alpha_i - \sum_{\i=1}^{3} \beta_i -2\gamma),\\
F_3(\mX,\mL) &= \frac{22848}{13047715}(\sum_{\i=1}^{3}\alpha_i - \sum_{\i=1}^{3} \beta_i -2\gamma),\\
F_4(\mX,\mL) &= \frac{5488}{2609543}(\sum_{\i=1}^{3}\alpha_i - \sum_{\i=1}^{3} \beta_i -2\gamma),\\
F_5(\mX,\mL) &= \frac{3360}{2609543}(\sum_{\i=1}^{3}\alpha_i - \sum_{\i=1}^{3} \beta_i -2\gamma),\\
F_6(\mX,\mL) &= \frac{3424}{13047715}(\sum_{\i=1}^{3}\alpha_i - \sum_{\i=1}^{3} \beta_i -2\gamma),\\
F_7(\mX,\mL) &= \frac{672}{13047715}(\sum_{\i=1}^{3}\alpha_i - \sum_{\i=1}^{3} \beta_i -2\gamma).
\end{align*}

\bibliographystyle{amsalpha}

  \end{document}